\documentclass[reqno]{amsart}
\usepackage[reset,a4paper,vmargin=5truecm,hmargin=4truecm]{geometry}
\usepackage{amsmath,amsthm,amssymb}

\usepackage{graphicx}
\usepackage{url}
\usepackage{cite}
\usepackage{hyperref}

\usepackage{fancyhdr}
\usepackage{lscape}

\newtheorem{thm}{Theorem}

\newtheorem{cor}{Corollary}
\newtheorem{prop}{Proposition}

\newtheorem{sta}{Statement}

\newtheorem{nameth}{Euler's Theorem for Tilings}

\begin{document}

\title[Properties of Strongly Balanced Tilings by Convex Polygons]{Properties of 
Strongly Balanced Tilings by Convex Polygons}

\address{The Interdisciplinary Institute of Science, Technology and Art\\
Suzukidaini-building 211, 2-5-28 Kitahara, Asaka-shi, Saitama, 351-0036, 
Japan}
\email{ismsugi@gmail.com}
\author{Teruhisa Sugimoto}
\maketitle	

\begin{abstract}
Every normal periodic tiling is a strongly balanced 
tiling. The properties of periodic tilings by convex polygons are rearranged 
from the knowledge of strongly balanced tilings. From the results, we show 
the properties of representative periodic tilings by a convex pentagonal 
tile.
\end{abstract}

\section{Introduction}
\label{section1}

A collection of sets (the `tiles') is a\textit{ tiling} (or tessellation) of the plane if 
their union covers the whole plane but the interiors of different tiles are 
disjoint. If all the tiles in a tiling are of the same size and shape, then 
the tiling is \textit{monohedral}. 
In this study, a polygon that admits a monohedral tiling is a \textit{polygonal tile}. 
A tiling by convex polygons is \textit{edge-to-edge} if any two convex polygons 
in a tiling are either disjoint or share one vertex (or an entire edge) in common. 
A tiling is periodic if it coincides with its translation by a nonzero vector. The 
unit that can generate a periodic tiling by translation only is known as the 
\textit{fundamental region}~\cite{G_and_S_1987, Sugimoto_2012a, Sugimoto_2015, 
Sugimoto_NoteTP, Sugimoto_APTCP}.

\begin{figure}[htbp]
 \centering\includegraphics[width=13cm,clip]{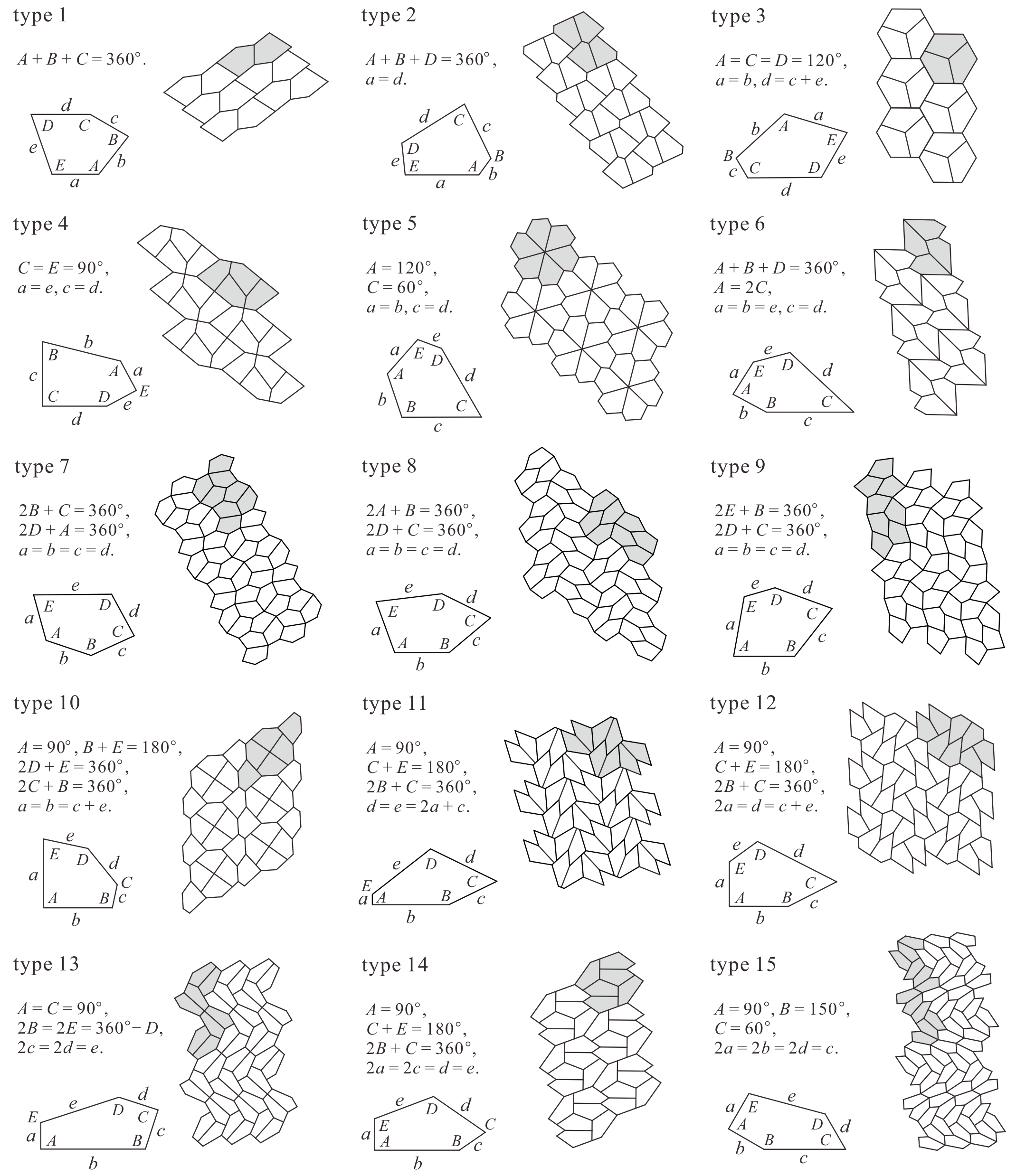} 
  \caption{{\small 
Fifteen types of convex pentagonal tiles. If a convex pentagon can generate 
a monohedral tiling and is not a new type, it belongs to at least one of types 1--15.
Each of the convex pentagonal tiles is defined by some conditions between 
the lengths of the edges and the magnitudes of the angles, but some degrees 
of freedom remain. For example, a convex pentagonal tile belonging to type 1 
satisfies that the sum of three consecutive angles is equal to $360^\circ$. 
This condition for type 1 is expressed as $A + B + C = 360^\circ$ in this figure. 
The pentagonal tiles of types 14 and 15 have one degree of freedom, that of 
size. For example, the value of $C$ of the pentagonal tile of type 14 is $\cos 
^{ - 1}((3\sqrt {57} - 17) / 16) \approx 1.2099 \; \mbox{rad} \approx 69.32^\circ$.
}
\label{fig1}
}
\end{figure}

In the classification problem of convex polygonal tiles, only the pentagonal 
case is open. At present, fifteen types of convex pentagon tiles are known 
(see Figure~\ref{fig1}) but it is not known whether this list is 
complete~\cite{Gardner_1975a, Gardner_1975b, G_and_S_1987, Hallard_1991,
Hirshh_1985, Kershner_1968, Klamkin_1980, Mann_2015, Reinhardt_1918, 
Schatt_1978, Stein_1985, Sugimoto_2012a, Sugimoto_2015, Sugi_Ogawa_2006}. 
However, it has been proved that a convex pentagonal tile that can generate 
an edge-to-edge tiling belongs to at least one of the eight known 
types~\cite{Bagina_2011, Sugimoto_2012b, Sugimoto_NoteTP, Sugimoto_2016}. 
We are interested in the problem of convex pentagonal tiling (i.e., the complete 
list of types of convex pentagonal tile, regardless of edge-to-edge and 
non-edge-to-edge tilings). However, the solution of the problem is not easy. 
Therefore, we will first treat only convex pentagonal tiles that admit at 
least one periodic tiling\footnote{ We know as a fact that the 15 types of 
convex pentagonal tile admit at least one periodic tiling. From this, we find 
that the convex pentagonal tiles that can generate an edge-to-edge tiling 
admit at least one periodic tiling~\cite{G_and_S_1987, Mann_2015,  Sugimoto_2012a, 
Sugimoto_2012b, Sugimoto_2015,  Sugimoto_NoteTP, Sugimoto_APTCP, Sugimoto_2016}. 
On the other hand, there is no proof that they admit at least one periodic 
tiling without using this fact. That is, there is no assurance yet that all 
convex pentagonal tiles admit at least one periodic tiling. In the solution 
of the problem of convex pentagonal tiling, it is necessary to consider whether 
there is a convex polygonal tile that admits infinitely many tilings of the 
plane, none of which is periodic.}. As such, we consider that the properties 
of the periodic tilings by convex polygons should be rearranged. From 
statement 3.4.8 (``every normal periodic tiling is strongly balanced'') in 
\cite{G_and_S_1987}, we see that periodic tilings by convex polygonal tiles 
(i.e., monohedral periodic tilings by convex polygons) are contained in 
the strongly balanced tilings. The definitions of normal and strongly balanced 
tilings are given in Section~\ref{section2}. In this paper, the properties of 
strongly balanced tilings by convex polygons are presented from the 
knowledge of strongly balanced tilings in general. That is, the properties 
correspond to those of periodic tilings by a convex polygonal tile.

\renewcommand{\figurename}{{\small Figure.}}

\renewcommand{\figurename}{{\small Figure.}}
\begin{figure}[htbp]
 \centering\includegraphics[width=11cm,clip]{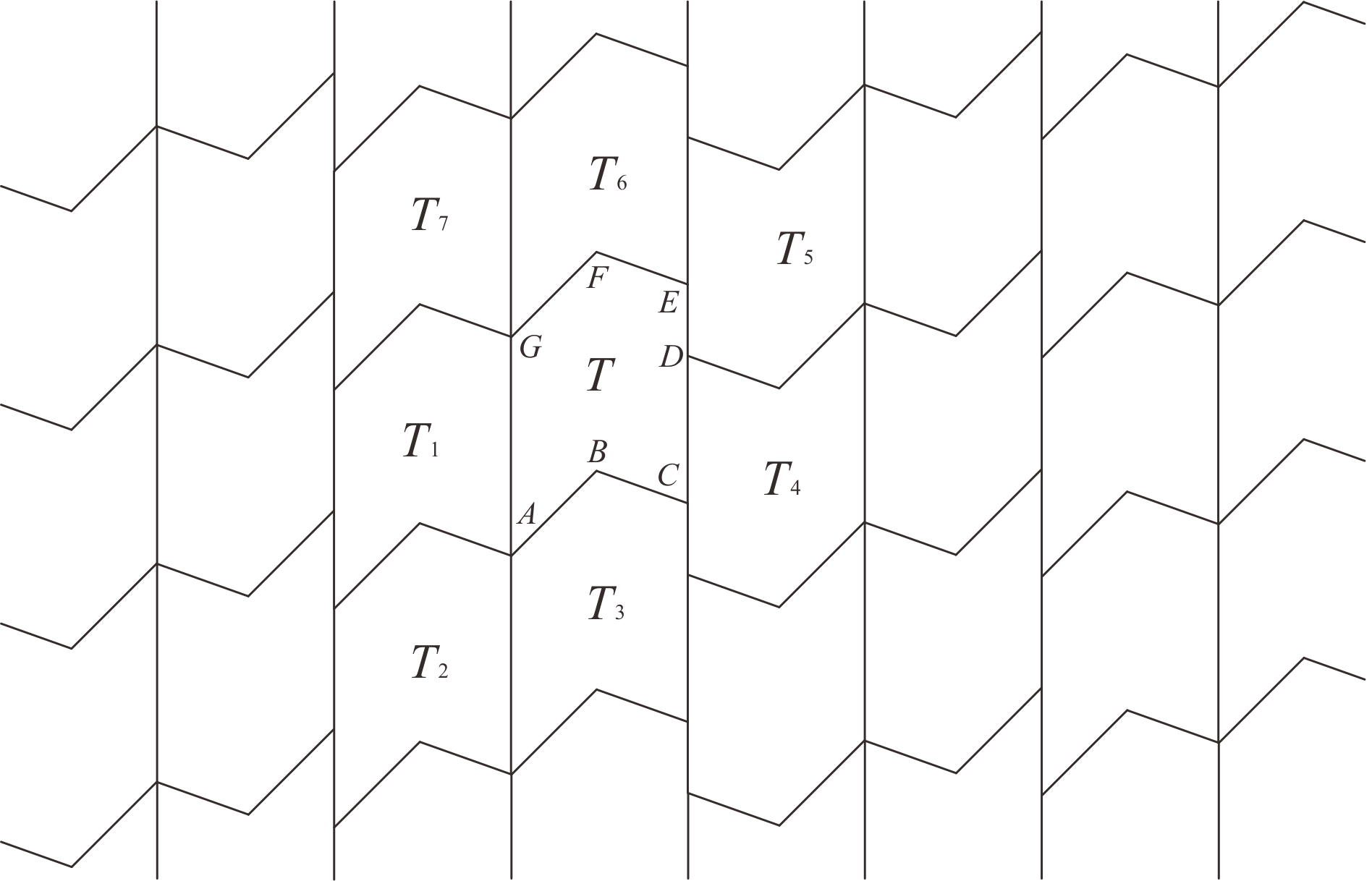} 
  \caption{{\small 
The differences between corners and vertices, sides and edges, 
adjacents, and neighbors. The points $A$, $B$, $C$, $E$, $F$, and $G$ are corners 
of the tile $T$; but $A$, $C$, $D$, $E$, and $G$ are vertices of the tiling (we note 
that the valence of vertices $A$ and $G$ is four, and the valence of vertices 
$C$, $D$, and $E$ is three). The line segments \textit{AB}, \textit{BC}, \textit{CE}, 
\textit{EF}, \textit{FG,} and \textit{GA} are sides of $T$, while \textit{AC}, 
\textit{CD}, \textit{DE}, \textit{EG}, and \textit{GA} are edges of the tiling. 
The tiles $T_{1}$, $T_{3}$, $T_{4}$, $T_{5}$, and $T_{6}$ are adjacents (and 
neighbors) of $T$, whereas tiles $T_{2}$ and $T_{7}$ are neighbors (but not 
adjacents) of $T$ \cite{G_and_S_1987}.
}
\label{fig2}
}
\end{figure}

\section{Preparation}
\label{section2}

Definitions and terms of this section quote from \cite{G_and_S_1987}. 

Terms ``vertices'' and ``edges'' are used by both polygons and tilings. In order 
not to cause confusion, \textit{corners} and \textit{sides} are referred to instead 
of vertices and the edges of polygons, respectively. At a vertex of a polygonal tiling, 
corners of two or more polygons meet and the number of polygons meeting at 
the vertex is called the \textit{valence} of the vertex , and is at least three (see 
Figure~\ref{fig2}). Therefore, an edge-to-edge tiling by polygons is such that the 
corners and sides of the polygons in a tiling coincide with the vertices and edges 
of the tiling.

Two tiles are called \textit{adjacent} if they have an edge in common, and then 
each is called an adjacent of the other. On the other hand, two tiles are called 
\textit{neighbors} if their intersection is nonempty (see Figure~\ref{fig2}).

There exist positive numbers $u$ and $U$ such that any tile contains a certain 
disk of radius $u$ and is contained in a certain disk of radius $U$ in which case 
we say the tiles in tiling are \textit{uniformly bounded}.

A tiling $\Im$ is called \textit{normal} if it satisfies following conditions: 
(i) every tiles of $\Im $ is a topological disk; (ii) the intersection of every two 
tiles of $\Im$ is a connected set, that is, it does not consist of two (or 
more) distinct and disjoint parts; (iii) the tiles of $\Im$ are uniformly 
bounded.

Let $D(r,M)$ be a closed circular disk of radius $r$, centered at any point 
$M$ of the plane. Let us place $D(r,M)$ on a tiling, and let $F_{1}$ and $F_{2}$ 
denote the set of tiles contained in $D(r,M)$ and the set of meeting 
boundary of $D(r,M)$ but not contained in $D(r,M)$, respectively. In 
addition, let $F_{3}$ denote the set of tiles surrounded by these in 
$F_{2}$ but not belonging to $F_{2}$. The set $F_1 \cup F_2 \cup F_3 $ of 
tiles is called the patch $A(r,M)$ of tiles generated by $D(r,M)$.

For a given tiling $\Im$, we denote by $v(r,M)$, $e(r,M)$, and $t(r,M)$ the 
numbers of vertices, edges, and tiles in $A(r,M)$, respectively. The tiling 
$\Im$ is called \textit{balanced} if it is normal and satisfies the following 
condition: the limits

\[
\mathop {\lim }\limits_{r \to \infty } \frac{v(r,M)}{t(r,M)}\quad 
\mbox{and}\quad \mathop {\lim }\limits_{r \to \infty } 
\frac{e(r,M)}{t(r,M)}
\]

\noindent
exist. Note that $v(r,M) - e(r,M) + t(r,M) = 1$ is called Euler's Theorem 
for Planar Maps.

\begin{sta}[Statement 3.3.13 in \cite{G_and_S_1987}] 
\label{sta1}
Every normal periodic tiling is balanced.
\end{sta}

\begin{nameth}[Statement 3.3.3 in \cite{G_and_S_1987}] 
\label{EulerThm}
For any normal tiling $\Im$, if one of the limits 
$v(\Im ) = \mathop {\lim }\limits_{r \to \infty } \frac{v(r,M)}{t(r,M)}$ or
$e(\Im ) = \mathop {\lim }\limits_{r \to \infty } \frac{e(r,M)}{t(r,M)}$ exists 
and is finite, then so does the other. Thus the tiling is balanced and, 
moreover,

\begin{equation}
\label{equqtion1}
v(\Im ) = e(\Im ) - 1.
\end{equation}

\end{nameth}

For a given tiling $\Im$, we write $t_h (r,M)$ for the number of tiles with 
$h$ adjacents in $A(r,M)$, and $v_j (r,M)$ for the numbers of $j$-valent 
vertices in $A(r,M)$. Then the tiling $\Im$ is called \textit{strongly balanced} 
if it is normal and satisfies the following condition: all the limits

\[
t_h (\Im ) = \mathop {\lim }\limits_{r \to \infty } \frac{t_h 
(r,M)}{t(r,M)}\quad \mbox{and}\quad v_j (\Im ) = \mathop {\lim }\limits_{r 
\to \infty } \frac{v_j (r,M)}{t(r,M)}
\]

\noindent
exist. Then,

\begin{equation}
\label{equqtion2}
\sum\limits_{h \ge 3} {t_h (\Im ) = 1} \quad \mbox{and}\quad v(\Im ) 
= \sum\limits_{j \ge 3} {v_j (\Im )}
\end{equation}

\noindent
hold. Therefore, every strongly balanced tiling is necessarily balanced. 

When $\Im$ is strongly balanced, we have

\begin{equation}
\label{equqtion3}
2e(\Im ) = \sum\limits_{j \ge 3} {j \cdot v_j (\Im )} 
= \sum\limits_{h \ge 3} {h \cdot t_h (\Im )} .
\end{equation}

In addition, as for strongly balanced tiling, following properties are 
known.

\begin{sta}[Statement 3.4.8 in \cite{G_and_S_1987}]
\label{sta2}
Every normal periodic tiling is strongly balanced.
\end{sta}

\begin{sta}[Statement 3.5.13 in \cite{G_and_S_1987}]
\label{sta3}
For each strongly balanced tiling $\Im $ we have

\begin{equation}
\label{equqtion4}
\frac{1}{\sum\limits_{j \ge 3} {j \cdot w_j (\Im )} } + 
\frac{1}{\sum\limits_{h \ge 3} {h \cdot t_h (\Im )} } = \frac{1}{2}
\end{equation}

\noindent
where

\[
w_j (\Im ) = \frac{v_j (\Im )}{v(\Im )}.
\]
\end{sta}

Thus $w_j (\Im )$ can be interpreted as that fraction of the total number of 
vertices in $\Im$ which have valence $j$, and $\sum\limits_{j \ge 3} {j \cdot 
w_j (\Im )} $ is the \textit{average} valence taken over all the vertices. Since 
$\sum\limits_{h \ge 3} {t_h (\Im ) = 1} $ there is a similar interpretation 
of $\sum\limits_{h \ge 3} {h \cdot t_h (\Im )} $: it is the \textit{average} 
number of adjacents of the tiles, taken over all the tiles in $\Im$. 
Since the valence of the vertex is at least three,

\begin{equation}
\label{equqtion5}
\sum\limits_{j \ge 3} {j \cdot w_j (\Im )} \ge 3.
\end{equation}

\begin{sta}[Statement 3.5.6 in \cite{G_and_S_1987}]
\label{sta4}
In every strongly balanced tiling $\Im$ we have

\[
2\sum\limits_{j \ge 3} {(j - 3) \cdot v_j (\Im )} + \sum\limits_{h \ge 3} 
{(h - 6) \cdot t_h (\Im ) = 0} ,
\]

\[
\sum\limits_{j \ge 3} {(j - 4) \cdot v_j (\Im )} + \sum\limits_{h \ge 3} {(h 
- 4) \cdot t_h (\Im ) = 0} ,
\]

\[
\sum\limits_{j \ge 3} {(j - 6) \cdot v_j (\Im )} + 2\sum\limits_{h \ge 3} 
{(h - 3) \cdot t_h (\Im ) = 0} .
\]

\end{sta}

\section{Consideration and Discussion}
\label{section3}

A polygon with $n$ sides and $n$ corners is referred to as an \textit{n-gon}. 
For the discussion below, note that $n$-gons in a strongly balanced tiling do 
not need to be congruent (i.e., the tiling does not need to be monohedral).

\subsection{Case of convex $n$-gons}
\label{subsection3_1}

Let $\Im _n^{sb} $ be a strongly balanced tiling by convex $n$-gons.

\begin{prop}\label{prop1}
$\sum\limits_{h \ge 3} {h \cdot t_h (\Im _n^{sb} )} \le 6$.
\end{prop}

\noindent
\textbf{\textit{Proof}.} 
From (\ref{equqtion4}), we have that

\[
\frac{2\sum\limits_{h \ge 3} {h \cdot t{ }_h(\Im _n^{sb} )} }{\sum\limits_{h 
\ge 3} {h \cdot t_h (\Im _n^{sb} )} - 2} 
= \sum\limits_{j \ge 3} {j \cdot w_j (\Im _n^{sb} )} .
\]

\noindent
Since the valence of the vertex is at least three, i.e., 
$\sum\limits_{j \ge 3} {j \cdot w_j (\Im _n^{sb} )} \ge 3$,

\[
\frac{2\sum\limits_{h \ge 3} {h \cdot t{ }_h(\Im _n^{sb} )} }{\sum\limits_{h 
\ge 3} {h \cdot t_h (\Im _n^{sb} )} - 2} \ge 3 \quad .
\]

\noindent
Therefore, we obtain Proposition~\ref{prop1}. 
\hspace{7cm} $\square$

\bigskip

From Proposition~\ref{prop1}, there is no strongly balanced tiling that is formed by 
convex $n$-gons for $n \ge 7$, since the average number of adjacents is greater 
than six. Note that the number of sides of all convex polygons in $\Im _n^{sb}$ 
does not have the same necessity. For example, there is no strongly balanced 
tiling by convex 6-gons and convex 8-gons, and there is no strongly balanced 
tiling by only convex 5-gons whose number of adjacents is seven or more.

Note that Proposition~\ref{prop1} is not a proof that there is no convex 
polygonal tile with seven or more sides. If there were a proof that all convex 
polygonal tiles admit at least one periodic tiling, it could be used to 
prove that there is no convex polygonal tile with seven or more sides from 
Proposition~\ref{prop1}.

\begin{prop}\label{prop2}
$3 \le \sum\limits_{j \ge 3} {j \cdot w_j (\Im _n^{sb} )} \le \frac{2n}{n - 2}$.
\end{prop}

\noindent
\textbf{\textit{Proof}.}  
From (\ref{equqtion4}) and $\sum\limits_{h \ge n} {h \cdot t_h (\Im _n^{sb} )} \ge n$,

\[
\frac{2\sum\limits_{j \ge 3} {j \cdot w{ }_j(\Im _n^{sb} )} }{\sum\limits_{j 
\ge 3} {j \cdot w_j (\Im _n^{sb} )} - 2} = \sum\limits_{h \ge n} {h \cdot 
t_h (\Im _n^{sb} )} \ge n.
\]

\noindent
Therefore, from the above inequality and (\ref{equqtion5}), we obtain 
Proposition~\ref{prop2}. 
\hspace{1.6cm} $\square$

\bigskip

Let $\Im _n^{sbe} $ be a strongly balanced edge-to-edge tiling by convex 
$n$-gons.

\begin{prop}\label{prop3}
$\sum\limits_{j \ge 3} {j \cdot w_j (\Im _n^{sbe} )} = \frac{2n}{n - 2}$.
\end{prop}

\noindent
\textbf{\textit{Proof}.}  
The number of adjacents of all convex $n$-gons in $\Im _n^{sb} $ is equal to $n$. 
That is, $\sum\limits_{h \ge n} {h \cdot t_h (\Im _n^{sbe} )} = n$. Then, 
(\ref{equqtion4}) is $\frac{1}{\sum\limits_{j \ge 3} {j \cdot w_j (\Im _n^{sbe} )} } 
+ \frac{1}{n} = \frac{1}{2}$. Therefore, we obtain Proposition~\ref{prop3}. 
\hspace{11.4cm} $\square$

\subsection{Case of convex hexagons}
\label{subsection3_2}

As for $\Im _6^{sb}$ (i.e., a strongly balanced tiling by convex hexagons 
(6-gon)), the number of adjacents of each convex hexagon should be greater 
than or equal to six (i.e., $h \ge 6)$. Therefore,

\begin{equation}
\label{equqtion6}
\sum\limits_{h \ge 6} {t_h (\Im _6^{sb} ) = t_6 (\Im _6^{sb} ) + 
\sum\limits_{h \ge 7} {t_h (\Im _6^{sb} )} } = 1.
\end{equation}

\noindent
On the other hand, from Proposition~\ref{prop1}, we have

\begin{equation}
\label{equqtion7}
\sum\limits_{h \ge 6} {h \cdot t_h (\Im _6^{sb} )} = 6 \cdot t_6 (\Im 
_6^{sb} ) + \sum\limits_{h \ge 7} {h \cdot t_h (\Im _6^{sb} )} \le 6.
\end{equation}

\begin{prop}\label{prop4}
$\sum\limits_{j \ge 3} {j \cdot w_j (\Im _6^{sb} )} = 3$.
\end{prop}

\noindent
\textbf{\textit{Proof}.} 
From (\ref{equqtion6}) and (\ref{equqtion7}),

\[
6\,\left( {1 - \sum\limits_{h \ge 7} {t_h (\Im _6^{sb} )} } \right) + 
\sum\limits_{h \ge 7} {h \cdot t_h (\Im _6^{sb} )} \le 6.
\]

\noindent
Hence, we obtain $\sum\limits_{h \ge 7} {(h - 6) \cdot t_h (\Im _6^{sb} )} \le 0$. 
On the other hand, $\sum\limits_{h \ge 7} {(h - 6) \cdot t_h (\Im _6^{sb} )} \ge 0$ 
holds because $t_h (\Im _6^{sb} ) = \mathop {\lim }\limits_{r \to \infty } 
\frac{t_h (r,M)}{t(r,M)} \ge 0$. From these inequalities,

\[
\sum\limits_{h \ge 7} {(h - 6) \cdot t_h (\Im _6^{sb} )} = 0.
\]

\noindent
Therefore, $t_h (\Im _6^{sb} ) = 0$ for $h \ne 6$. That is, we have 

\[
t_6 (\Im _6^{sb} ) = 1,
\quad
\sum\limits_{h \ge 6} {h \cdot t_h (\Im _6^{sb} )} = 6
\quad \mbox{and}\quad
\sum\limits_{h 
\ge 7} {t_h (\Im _6^{sb} )} = \sum\limits_{h \ge 7} {h \cdot t_h (\Im _6^{sb} )} = 0.
\]

\noindent
Thus, from these relationships and (\ref{equqtion4}), we arrive at 
Proposition~\ref{prop4}. 
\hspace{2.4cm} $\square$

\bigskip

From Statement~\ref{sta2}, a monohedral periodic tiling by a convex hexagon 
is strongly balanced. If a fundamental region in a monohedral periodic tiling 
by a convex hexagon has vertices with valences of four or more, 
$\sum\limits_{j \ge 3} {j \cdot w_j (\Im _6^{sb} )} > 3$ and 
$\sum\limits_{h \ge 6} {h \cdot t_h (\Im _6^{sb} )} < 6$ from (\ref{equqtion4}), 
which is a contradiction of Proposition~\ref{prop4}. Thus, we have the 
following corollary.

\begin{cor}\label{cor1}
A monohedral periodic tiling by a convex hexagon is an edge-to-edge tiling 
with only $3$-valent vertices.
\end{cor}

\begin{figure}[htbp]
 \centering\includegraphics[width=13cm,clip]{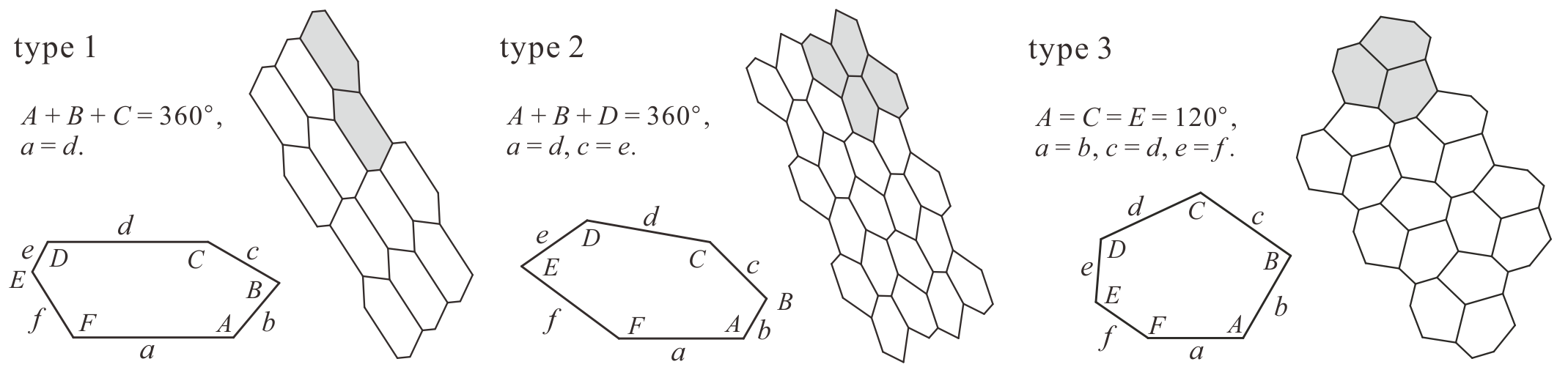} 
  \caption{{\small 
Three types of convex hexagonal tiles. If a convex hexagon can 
generate a monohedral tiling, it belongs to at least one of types 1--3. 
The pale gray hexagons in each tiling indicate the fundamental region.
}
\label{fig3}
}
\end{figure}

\begin{figure}[htbp]
 \centering\includegraphics[width=10cm,clip]{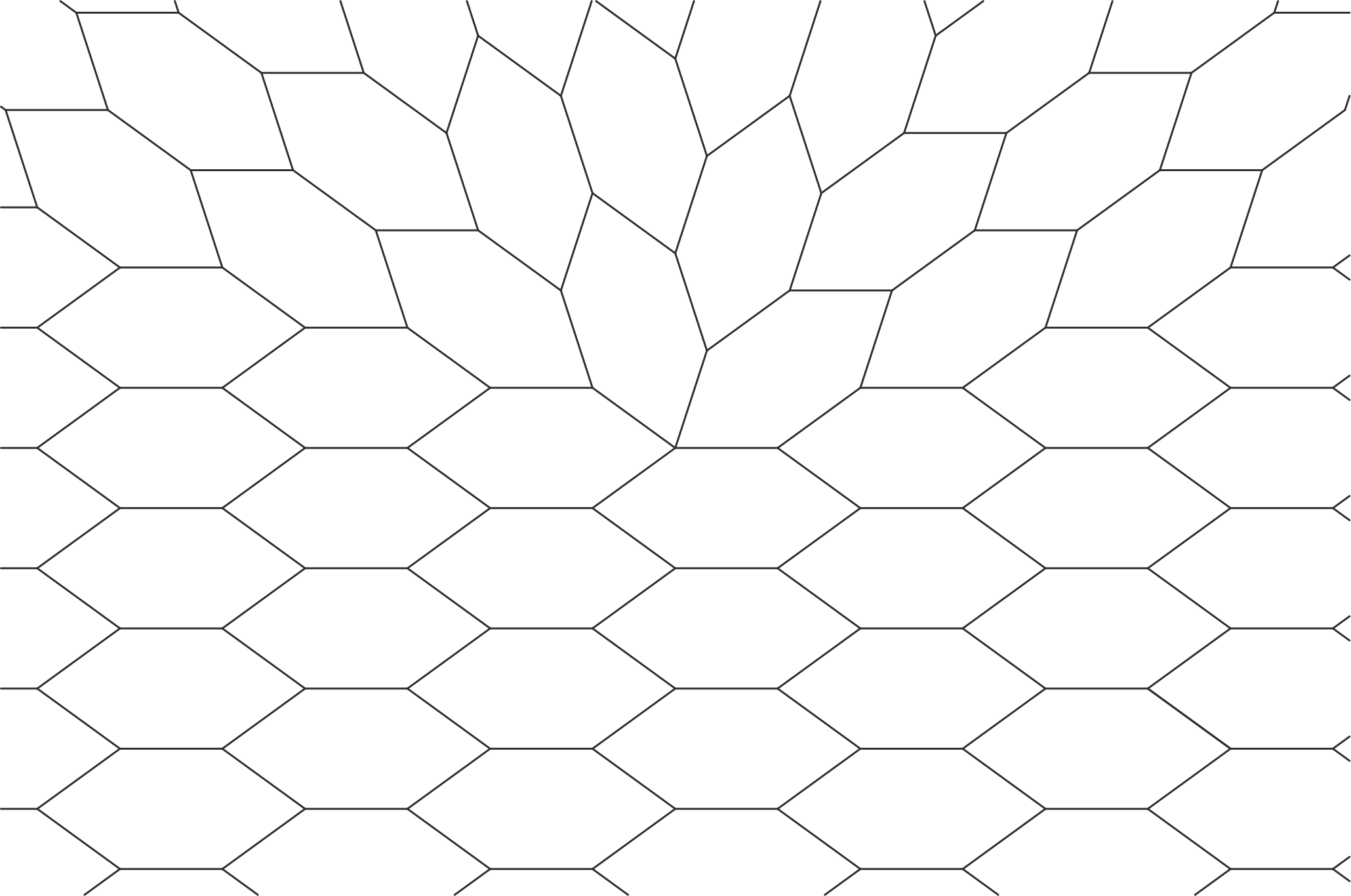} 
  \caption{{\small 
Monohedral tilings by convex hexagons with vertices of valence equal 
to four.}
\label{fig4}
}
\end{figure}

It is well known that convex hexagonal tiles (i.e., convex hexagons that 
admit a monohedral tiling) belong to at least one of the three types shown 
in Figure~\ref{fig3}. That is, convex hexagonal tiles admit at least one periodic 
edge-to-edge tiling whose valence is three at all vertices. In fact, the 
representative tilings of the three types in Figure~\ref{fig3} are periodic 
edge-to-edge tilings whose valence is three at all vertices. From Corollary 
1 and the fact that the valence of vertices is at least three, it might be 
considered that the valence of all vertices in monohedral tilings by convex 
hexagons is three; however, that is not true. For example, as shown Figure~\ref{fig4}, 
there are monohedral tilings by convex hexagons with vertices of valence 
equal to four (note that a monohedral tiling is not always a periodic 
tiling). However, it is clear that the convex hexagonal tiles of Figure~\ref{fig4} can 
generate a periodic edge-to-edge tiling in which the valence of all vertices 
is three.

\renewcommand{\figurename}{{\small Figure.}}

\renewcommand{\figurename}{{\small Figure.}}

\subsection{Case of convex pentagons}
\label{subsection3_3}

As for $\Im _5^{sb}$ (i.e., a strongly balanced tiling by convex pentagons 
(5-gons)), the number of adjacents of each convex pentagon should be greater 
than or equal to five. From Proposition~\ref{prop1}, we obtain the following.

\begin{prop}\label{prop5}
$5 \le \sum\limits_{h \ge 5} {h \cdot t_h (\Im _5^{sb} )} \le 6$. 
\end{prop}

\noindent
Since $\sum\limits_{h \ge 5} {h \cdot t_h (\Im _5^{sb} )}$ is the average 
number of adjacents of the convex pentagons in $\Im _5^{sb}$, we obtain 
the following theorem~\cite{Sugimoto_2014}.

\begin{thm}\label{thm1}
A tiling  $\Im _5^{sb}$ contains a convex pentagon whose number of 
adjacents is five or six.
\end{thm}

\noindent
Then, we obtain the following propositions.

\begin{prop}\label{prop6}
$\frac{5}{2} \le e(\Im _5^{sb} ) \le 3$. 
\end{prop}

\begin{prop}\label{prop7}
$\frac{3}{2} \le v(\Im _5^{sb} ) \le 2$.
\end{prop}

\noindent
\textbf{\textit{Proof of Propositions $\ref{prop6}$ and $\ref{prop7}$}.} 
From (\ref{equqtion3}) and Proposition~\ref{prop5}, we have that 
$5 \le 2e(\Im _5^{sb} ) \le 6$. Since each strongly balanced tiling is 
necessarily balanced, $v(\Im _5^{sb} ) = e(\Im _5^{sb} ) - 1$ holds from 
Euler's Theorem for Tilings  (see P.\pageref{EulerThm}). Therefore, we 
have that $5 \le 2(v(\Im _5^{sb} ) +1) \le 6$. Thus, we obtain
Propositions $\ref{prop6}$ and $\ref{prop7}$.\hspace{3.2cm} $\square$

\bigskip

\begin{prop}\label{prop8}
$3 \le \sum\limits_{j \ge 3} {j \cdot w_j (\Im _5^{sb} )} \le \frac{10}{3}$.
\end{prop}

\noindent
\textbf{\textit{Proof}.} 
It is clear from Proposition~\ref{prop2}. 
\hspace{6.6cm} $\square$

\bigskip

\begin{prop}\label{prop9}
$0 \le \sum\limits_{h \ge 7} {(h - 6) \cdot t_h (\Im _5^{sb} )} \le t_5 (\Im _5^{sb} ) \le 1$.
\end{prop}

\noindent
\textbf{\textit{Proof}.} 
From Proposition~\ref{prop5},

\begin{equation}
\label{equqtion8}
5 \le \sum\limits_{h \ge 5} {h \cdot t_h (\Im _5^{sb} )} = 5t_5 (\Im _5^{sb} 
) + 6t_6 (\Im _5^{sb} ) + \sum\limits_{h \ge 7} {h \cdot t_h (\Im _5^{sb} )} 
\le 6.
\end{equation}

\noindent
From (\ref{equqtion2}),

\[
\sum\limits_{h \ge 3} {t_h (\Im _5^{sb} )} = t_5 (\Im _5^{sb} ) + t_6 (\Im _5^{sb} ) 
+ \sum\limits_{h \ge 7} {t_h (\Im _5^{sb} )} = 1,
\]

\begin{equation}
\label{equqtion9}
t_6 (\Im _5^{sb} ) = 1 - t_5 (\Im _5^{sb} ) - \sum\limits_{h \ge 7} {t_h (\Im _5^{sb} )} .
\end{equation}

\noindent
From (\ref{equqtion8}), (\ref{equqtion9}) and $0 \le t_h (\Im _5^{sb} ) \le 1$, we obtain 
the inequality of Proposition~\ref{prop9}. 
\hspace{0.6cm} $\square$

\bigskip

From Proposition~\ref{prop9}, we obtain the following theorem~\cite{Sugimoto_2014}.

\begin{thm}\label{thm2}
A tiling $\Im _5^{sb}$ that satisfies $\sum\limits_{h \ge 7} {t_h (\Im _5^{sb} )} > 0$ 
contains a convex pentagon whose number of adjacents is five.
\end{thm}

As for $\Im _5^{sbe} $ (i.e., a strongly balanced edge-to-edge tiling by 
convex pentagons), we have the following proposition.

\begin{prop}\label{prop10}
$t_5 (\Im _5^{sbe} ) = 1$, $v(\Im _5^{sbe} ) = \frac{3}{2}$, 
$e(\Im _5^{sbe} ) = \frac{5}{2}$, and 
$\sum\limits_{j \ge 3} {j \cdot w_j (\Im _5^{sbe} )} = \frac{10}{3}$.
\end{prop}

\noindent
\textbf{\textit{Proof}.} 
It is clear from Proposition~\ref{prop3}, Euler's Theorem for Tilings  
(see P.\pageref{EulerThm}), and $t_h (\Im _5^{sbe} ) = 0$ for $h \ne 5$. 
\hspace{9cm} $\square$

\bigskip

Next, we have other propositions as follows.

\begin{prop}\label{prop11}
$v_3 (\Im _5^{sb} ) = 2 + \sum\limits_{j \ge 4} {(2 - j) \cdot v_j (\Im _5^{sb} )} $.
\end{prop}

\noindent
\textbf{\textit{Proof}.} 
From (\ref{equqtion1}) and the definition of $v(\Im )$, we have

\begin{equation}
\label{equqtion10}
e(\Im _5^{sb} ) = v(\Im _5^{sb} ) + 1 
= \sum\limits_{j \ge 3} {v_j (\Im _5^{sb} ) + 1 } 
=  v_3 (\Im _5^{sb} ) + \sum\limits_{j \ge 4} {v_j (\Im _5^{sb} )} + 1.
\end{equation}

\noindent
On the other hand, from (\ref{equqtion3}),

\begin{equation}
\label{equqtion11}
2e(\Im _5^{sb} ) = \sum\limits_{j \ge 3} {j \cdot v_j (\Im _5^{sb} ) = 
3v_3 (\Im _5^{sb} ) + \sum\limits_{j \ge 4} {j \cdot v_j (\Im _5^{sb} )} } 
\end{equation}

\noindent
holds. Therefore, from (\ref{equqtion10}) and (\ref{equqtion11}), we have

\[
2\,\left( {v_3 (\Im _5^{sb} ) + \sum\limits_{j \ge 4} {v_j (\Im _5^{sb} ) + 
1} } \right) = 3v_3 (\Im _5^{sb} ) + \sum\limits_{j \ge 4} {j \cdot v_j (\Im 
_5^{sb} )} .
\]

\noindent
Thus, we obtain Proposition~\ref{prop11}. 
\hspace{7.7cm} $\square$

\bigskip

\begin{prop}\label{prop12}
$v_3 (\Im _5^{sbe} ) = \sum\limits_{j \ge 4} {(3j - 10)} \cdot v_j (\Im _5^{sbe} )$.
\end{prop}

\noindent
\textbf{\textit{Proof}.} 
From Proposition~\ref{prop10} and the definitions of 
$w_j \left( \Im \right)$ and $v\left( \Im \right)$,

\[
\sum\limits_{j \ge 3} {j \cdot w_j (\Im _5^{sbe} ) = } \frac{\sum\limits_{j 
\ge 3} {j \cdot v_j (\Im _5^{sbe} )} }{v(\Im _5^{sbe} )} = \frac{3v_3 (\Im 
_5^{sbe} ) + \sum\limits_{j \ge 4} {j \cdot v_j (\Im _5^{sbe} )} }{v_3 (\Im 
_5^{sbe} ) + \sum\limits_{j \ge 4} {v_j (\Im _5^{sbe} )} } = \frac{10}{3}.
\]

\noindent
Thus, we obtain Proposition~\ref{prop12}. 
\hspace{7.7cm} $\square$

\bigskip

As for $v(\Im _5^{sb} )$, we have the following propositions.

\begin{prop}\label{prop13}
$v(\Im _5^{sb} ) = \sum\limits_{j \ge 3} {v_j (\Im _5^{sb} )} = 
\frac{1}{2}\sum\limits_{h \ge 5} {\left( {h - 2} \right)} \cdot t_h (\Im _5^{sb} )$.
\end{prop}

\begin{prop}\label{prop14}
$v(\Im _5^{sb} ) = \frac{1}{2} + 
\frac{1}{2}\sum\limits_{h \ge 5} {\left( {h - 3} \right)} \cdot t_h (\Im _5^{sb} )$.
\end{prop}

\begin{prop}\label{prop15}
$v(\Im _5^{sb} ) = 2 + \frac{1}{2}\sum\limits_{h \ge 5} {(h - 6) \cdot t_h (\Im _5^{sb} )}$.
\end{prop}

\begin{prop}\label{prop16}
$v(\Im _5^{sb} ) = 2 - \sum\limits_{j \ge 4} {(j - 3) \cdot v_j (\Im _5^{sb} )} $.
\end{prop}

\noindent
\textbf{\textit{Proof of  Propositions $\ref{prop13}, \ref{prop14}, \ref{prop15},$ and 
$\ref{prop16}$}.} From the first equation in Statement~\ref{sta4},

\[
2\sum\limits_{j \ge 4} {(j - 3) \cdot v_j (\Im _5^{sb} )} + 
\sum\limits_{h \ge 5} {(h - 6) \cdot t_h (\Im _5^{sb} )} = 0.
\]

\noindent
Note that $t_3 (\Im _5^{sb} ) = t_4 (\Im _5^{sb} ) = 0$, since $\Im _5^{sb}$ 
is $h \ge 5$. The above equation is rearranged as

\begin{equation}
\label{equqtion12}
\sum\limits_{j \ge 4} {(j - 3) \cdot v_j (\Im _5^{sb} )} 
= - \frac{1}{2}\sum\limits_{h \ge 5} {(h - 6) \cdot t_h (\Im _5^{sb} )}.
\end{equation}

\noindent
From the second equation in Statement~\ref{sta4},

\[
 - v_3 (\Im _5^{sb} ) + \sum\limits_{j \ge 5} {(j - 4) \cdot v_j (\Im _5^{sb} )} 
+ \sum\limits_{h \ge 5} {(h - 4) \cdot t_h (\Im _5^{sb} )} = 0.
\]

\noindent
The above equation is rearranged as

\begin{equation}
\label{equqtion13}
\sum\limits_{j \ge 5} {(j - 4) \cdot v_j (\Im _5^{sb} )} = 
v_3 (\Im _5^{sb} ) - \sum\limits_{h \ge 5} {(h - 4) \cdot t_h (\Im _5^{sb} )} .
\end{equation}

\noindent
Then,

\begin{equation}
\label{equqtion14}
\begin{array}{l}
 \sum\limits_{j \ge 4} {(j - 3) \cdot v_j (\Im _5^{sb} )} = v_4 (\Im _5^{sb} ) 
+ \sum\limits_{j \ge 5} {(j - 4 + 1) \cdot v(\Im _5^{sb} )} \\ 
 \quad \quad \quad \quad \quad \quad \quad \quad \:\:\: = v_4 (\Im _5^{sb} ) + 
\sum\limits_{j \ge 5} {(j - 4) \cdot v(\Im _5^{sb} ) + \sum\limits_{j \ge 5} 
{v_j (\Im _5^{sb} )} }. \\ 
 \end{array}
\end{equation}

\noindent
By replacing $\sum\limits_{j \ge 4} {(j - 3) \cdot v_j (\Im _5^{sb} )} $ of 
(\ref{equqtion12}) and $\sum\limits_{j \ge 5} {(j - 4) \cdot v_j (\Im _5^{sb} )} $ 
of (\ref{equqtion13}) in (\ref{equqtion14}), the latter becomes

\[
 - \frac{1}{2}\sum\limits_{h \ge 5} {(h - 6) \cdot t_h (\Im _5^{sb} )} = 
v_4 (\Im _5^{sb} ) + v_3 (\Im _5^{sb} ) - \sum\limits_{h \ge 5} {(h - 4) \cdot 
t_h (\Im _5^{sb} )} + \sum\limits_{j \ge 5} {v_j (\Im _5^{sb} )} .
\]

\noindent
Simplifying both sides,

\[
\sum\limits_{j \ge 3} {v_j (\Im _5^{sb} ) = } \frac{1}{2}\sum\limits_{h \ge 5} 
{\left( { - h + 6 + 2h - 8} \right) \cdot t_h (\Im _5^{sb} ) = } 
\frac{1}{2}\sum\limits_{h \ge 5} {\left( {h - 2} \right) \cdot t_h (\Im _5^{sb} )} .
\]

\noindent
Thus, we obtain Proposition~\ref{prop13}. 

Next, from $\sum\limits_{j \ge 4} {(j - 3) \cdot v_j (\Im _5^{sb} )} $ and 
Proposition~\ref{prop11},

\[
\begin{array}{l}
 \sum\limits_{j \ge 4} {(j - 3) \cdot v_j (\Im _5^{sb} )} = - \sum\limits_{j \ge 4} {(2 - j) 
\cdot v_j (\Im _5^{sb} )} - \sum\limits_{j \ge 4} {v_j (\Im_5^{sb} )} \\ 
\quad \quad \quad \quad \quad \quad \quad \quad \:\:\:
 = 2 - v_3 (\Im _5^{sb} ) - \sum\limits_{j \ge 4} {v_j (\Im _5^{sb} )} \\ 
\quad \quad \quad \quad \quad \quad \quad \quad \:\:\:
= 2 - \sum\limits_{j \ge 3} {v_j (\Im _5^{sb} )}. \\ 
\end{array}
\]

\noindent
Thus, we obtain Proposition~\ref{prop16}. 

From (\ref{equqtion12}) and Proposition~\ref{prop16},

\[
2 - \sum\limits_{j \ge 3} {v_j (\Im _5^{sb} )} = 
- \frac{1}{2}\sum\limits_{h \ge 5} {(h - 6) \cdot t_h (\Im _5^{sb} )} .
\]

\noindent
Thus, we obtain Proposition~\ref{prop15}.

From the third equation in Statement~\ref{sta4},

\[
 - 3v_3 (\Im _5^{sb} ) + \sum\limits_{j \ge 4} {(j - 6) \cdot v_j (\Im _5^{sb} )} 
+ 2\sum\limits_{h \ge 5} {(h - 3) \cdot t_h (\Im _5^{sb} ) = 0}.
\]

\noindent
By replacing $v_3 (\Im _5^{sb} )$ of Proposition~\ref{prop11} in the above equation, 
it becomes

\begin{equation}
\label{equqtion15}
4\sum\limits_{j \ge 4} {(j - 3) \cdot v_j (\Im _5^{sb} )} = 
6 - 2\sum\limits_{h \ge 5} {(h - 3) \cdot t_h (\Im _5^{sb} )}.
\end{equation}

\noindent
Form Proposition~\ref{prop11} and (\ref{equqtion15}),

\[
\begin{array}{l}
v_3 (\Im _5^{sb} ) = 2 - \sum\limits_{j \ge 4} {v_j (\Im _5^{sb} )} - 
\sum\limits_{j \ge 4} {(j - 3) \cdot v_j (\Im _5^{sb} )}\\ \quad \quad \quad \;\; 
= \frac{1}{2} - \sum\limits_{j \ge 4} {v_j (\Im _5^{sb} )} 
+ \frac{1}{2}\sum\limits_{h \ge 5} {(h - 3) \cdot t_h (\Im _5^{sb} )}.
\end{array}
\]

\noindent
Simplifying,

\[
v_3 (\Im _5^{sb} ) + \sum\limits_{j \ge 4} {v_j (\Im _5^{sb} )} = 
\frac{1}{2} + \frac{1}{2}\sum\limits_{h \ge 5} {(h - 3) \cdot t_h (\Im 
_5^{sb} )} .
\]

\noindent
Thus, we obtain Proposition~\ref{prop14}. 
\hspace{7.7cm} $\square$
\bigskip

Here, we consider the case of $v(\Im _5^{sb} ) = \frac{3}{2}$ (i.e., the 
minimum case of $v(\Im _5^{sb} ))$. From Proposition~\ref{prop15}, we have

\[
2 + \frac{1}{2}\sum\limits_{h \ge 5} {(h - 6) \cdot t_h (\Im _5^{sb} )} 
= 2 - \frac{1}{2}t_5 ( {\Im _5^{sb} } ) + 
\frac{1}{2}\sum\limits_{h \ge 6} {(h - 6) \cdot t_h (\Im _5^{sb} )} = \frac{3}{2}.
\]

\noindent
Simplifying this equation,

\begin{equation}
\label{equqtion16}
t_5 ( {\Im _5^{sb} } ) = 1 + \sum\limits_{h \ge 6} {(h - 6) \cdot 
t_h (\Im _5^{sb} )} .
\end{equation}

\noindent
Since $\sum\limits_{h \ge 5} {t_h (\Im _5^{sb} ) = t_5 \left( {\Im _5^{sb} } \right) + } 
\sum\limits_{h \ge 6} {t_h (\Im _5^{sb} ) = 1} $, 
$\sum\limits_{h \ge 6} {(h - 6) \cdot t_h (\Im _5^{sb} )} $ in (\ref{equqtion16}) is equal 
to zero. That is, in the case of $v(\Im _5^{sb} ) = \frac{3}{2}$, 
$\sum\limits_{h \ge 6} {t_h (\Im _5^{sb} )} = 0$. Therefore, $v(\Im _5^{sb} 
) = \frac{3}{2}$ if and only if $\Im _5^{sb} $ is $\Im _5^{sbe} $.

\subsection{Properties of representative periodic tilings by a 
convex pentagonal tile}
\label{subsection3_4}

\begin{figure}[htbp]
 \centering\includegraphics[width=13cm,clip]{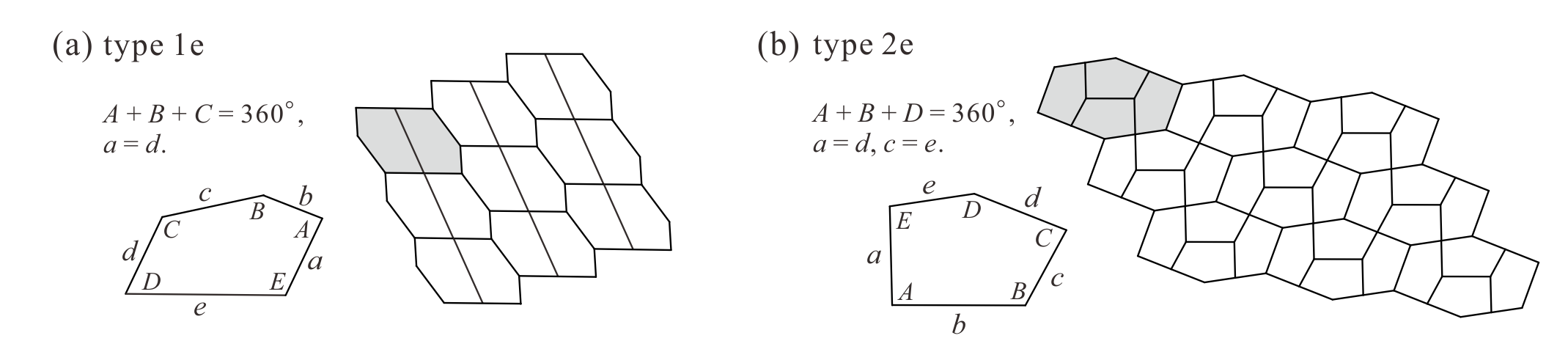} 
  \caption{{\small 
Examples of edge-to-edge tilings by convex pentagonal tiles that 
belong to types 1 or 2. The pale gray pentagons in each tiling indicate the 
fundamental region.}
\label{fig5}
}
\end{figure}

Let $\Im _5^{r(x)} $ be a representative periodic tiling by a convex 
pentagonal tile of type $x$. That is, $\Im _5^{r(x)} $ for $x = 1,\,\ldots ,\,15$ 
is a strongly balanced tiling by a convex pentagonal tile. 
Representative tilings of types 1 or 2 are generally non-edge-to-edge, as 
shown in Figure~\ref{fig1}. However, in special cases, the convex pentagonal tiles of 
types 1 or 2 can generate edge-to-edge tilings, as shown in Figure~\ref{fig5}. Here, 
the convex pentagonal tiles of (a) and (b) in Figure~\ref{fig5} are referred to as 
those of types 1e and 2e, respectively. Then, let $\Im _5^{r(1e)} $ and 
$\Im _5^{r(2e)} $ be representative edge-to-edge periodic tilings by a convex 
pentagonal tile of types 1e and 2e, respectively.

The properties of each tiling $\Im _5^{r(x)} $ are obtained from the results 
of Section~\ref{subsection3_3}, etc. Table 1 summarizes the results~\cite{Sugimoto_2014}. 
We can check that the representative periodic tilings of each type of convex 
pentagonal tile that can generate an edge-to-edge tiling satisfy 
Proposition~\ref{prop10}.

\renewcommand{\figurename}{{\small Figure.}}

\section{Conclusions}

In this paper, although it is accepted as a fact, it is proved that the only 
convex polygonal tiles that admit at least one periodic tiling are 
triangles, quadrangles, pentagons, and hexagons. As for the fact (proof) 
that the only convex polygonal tiles are triangles, quadrangles, pentagons, 
and hexagons, note that it is necessary to take except periodic tilings also 
into consideration. 

Although a solution to the problem of classifying the types of convex 
pentagonal tile is not yet in sight, we suggest that the properties that 
could lead to such a solution are those that are shown in Section~\ref{section3}.

\begin{landscape}

\begin{table}[htbp]

{\small 
\caption[Table 1.]{Properties of $\Im _5^{r(x)} $}
\label{table1}
}
\

\renewcommand{\arraystretch}{2.1}
\begin{tabular}
{p{90pt}|p{105pt}p{90pt}p{50pt}p{90pt}p{75pt}}

\hline
\hfil \hfil \raisebox{0ex}[0cm][0cm]{$x$ of $\Im _5^{r(x)} $}
& 
\hfil \hfil \raisebox{0.5ex}[0cm][0cm]{$t_h \Bigl(\Im _5^{r(x)} \Bigr)$}
& 
\hfil \hfil \raisebox{0.5ex}[0cm][0cm]{$v_j \Big(\Im _5^{r(x)} \Bigr)$}
& 
\hfil \hfil \raisebox{0.5ex}[0cm][0cm]{$2e \Big(\Im _5^{r(x)} \Bigr)$}
& 
\hfil \raisebox{0.5ex}[0cm][0cm]{$w_j \Big(\Im _5^{r(x)} \Bigr)$}
& 
\hfil \raisebox{1ex}[0cm][0cm]{$\sum\limits_{j \ge 3} {j \cdot w_j \Big(\Im _5^{r(x)} \Bigr)}  $} \\

\hline
\hfil \raisebox{-1.50ex}[0cm][0cm]{1e, 2e, 4,  6, 7, 8, 9}
& 
$t_5 = 1,$ \par $t_h = 0$ for $h \ne 5$
& 
$v_3 = 1, v_4 = \frac{1}{2},$ \par $v_j = 0$ for $j \ne 3\;\mbox{or}\;4$
& 
\hfil \raisebox{-1.50ex}[0cm][0cm]{5}
& 
$w_3 = \frac{2}{3}, w_4 = \frac{1}{3},$ \par $w_j = 0$ for $j \ne 3\;\mbox{or}\;4$
& 
\hfil \raisebox{-1.50ex}[0cm][0cm]{$\frac{10}{3} = 3.\dot {3}$} \\

\hline 
\hfil \raisebox{-1.50ex}[0cm][0cm]{5}
& 
$t_5 = 1,$ \par $t_h = 0$ for $h \ne 5$
& 
$v_3 = \frac{4}{3}, v_6 = \frac{1}{6},$ \par $v_j = 0$ for $j \ne 3\;\mbox{or}\;6$
&
\hfil \raisebox{-1.50ex}[0cm][0cm]{5}
& 
$w_3 = \frac{8}{9}, w_6 = \frac{1}{9},$ \par $w_j = 0$ for $j \ne 3\;\mbox{or}\;6$
& 
\hfil \raisebox{-1.50ex}[0cm][0cm]{$\frac{10}{3} = 3.\dot {3}$} \\

\hline
\hfil \raisebox{-1.50ex}[0cm][0cm]{1, 2, 3, 12}
& 
$t_6 = 1,$ \par $t_h = 0$ for $h \ne 6$
& 
$v_3 = 2,$ \par $v_j = 0$ for $j \ne 3$
& 
\hfil \raisebox{-1.50ex}[0cm][0cm]{6}
& 
$w_3 = 1,$ \par $w_j = 0$ for $j \ne 3$& 
\hfil \raisebox{-1.50ex}[0cm][0cm]{3} \\

\hline
\hfil \raisebox{-1.50ex}[0cm][0cm]{10}
& 
$t_5 = \frac{2}{3},t_7 = \frac{1}{3},$ \par $t_h = 0$ for $h \ne 5\;\mbox{or}\;7$
& 
$v_3 = \frac{5}{3},v_4 = \frac{1}{6},$ \par $v_j = 0$ for $j \ne 3\;\mbox{or}\;4$
& 
\hfil \raisebox{-1.50ex}[0cm][0cm]{$\frac{17}{3} = 5.\dot {6}$}
& 
$w_3 = \frac{10}{11},w_4 = \frac{1}{11},$ \par $w_j = 0$ for $j \ne 3\;\mbox{or}\;4$
& 
\hfil \raisebox{-1.50ex}[0cm][0cm]{$\frac{34}{11} = 3.\dot {0}\dot {9}$} \\

\hline
\hfil \raisebox{-1.50ex}[0cm][0cm]{11}
& 
$t_5 = t_7 = \frac{1}{2},$ \par $t_h = 0$ for $h \ne 5\;\mbox{or}\;7$
& 
$v_3 = 2,$ \par $v_j = 0$ for $j \ne 3$
& 
\hfil \raisebox{-1.50ex}[0cm][0cm]{6}
& 
$w_3 = 1,$ \par $w_j = 0$ for $j \ne 3$
& 
\hfil \raisebox{-1.50ex}[0cm][0cm]{3} \\

\hline
\hfil \raisebox{-1.50ex}[0cm][0cm]{13}
& 
$t_5 = t_6 = \frac{1}{2},$ \par $t_h = 0$ for $h \ne 5\;\mbox{or}\;6$
& 
$v_3 = \frac{3}{2},v_4 = \frac{1}{4},$ \par $v_j = 0$ for $j \ne 3\;\mbox{or}\;4$
& 
\hfil \raisebox{-1.50ex}[0cm][0cm]{$\frac{11}{2} = 5.5$}
& 
$w_3 = \frac{6}{7},w_4 = \frac{1}{7},$ \par $w_j = 0$ for $j \ne 3\;\mbox{or}\;4$
& 
\hfil \raisebox{-1.50ex}[0cm][0cm]{$\frac{22}{7} \approx 3.142...$} \\

\hline
\hfil \raisebox{-1.50ex}[0cm][0cm]{14}
& 
$t_5 = t_6 = t_7 = \frac{1}{3},$ \par $t_h = 0$ for $h \ne 5, 6, \;\mbox{or}\;7$
& 
$v_3 = 2$ \par $v_j = 0$ for $j \ne 3$
& 
\hfil \raisebox{-1.50ex}[0cm][0cm]{6}
& 
$w_3 = 1,$ \par $w_j = 0$ for $j \ne 3$
& 
\hfil \raisebox{-1.50ex}[0cm][0cm]{3} \\

\hline
\hfil \raisebox{-1.50ex}[0cm][0cm]{15}
& 
$t_5 = \frac{2}{3},t_6 = \frac{1}{3},$ \par $t_h = 0$ for $h \ne 5\;\mbox{or}\;6$
& 
$v_3 = \frac{4}{3},v_4 = \frac{1}{3},$ \par $v_j = 0$ for $j \ne 3\;\mbox{or}\;4$
& 
\hfil \raisebox{-1.50ex}[0cm][0cm]{$\frac{16}{3} = 5.\dot {3}$}
& 
$w_3 = \frac{4}{5},w_4 = \frac{1}{5},$ \par $w_j = 0$ for $j \ne 3\;\mbox{or}\;4$
& 
\hfil \raisebox{-1.50ex}[0cm][0cm]{$\frac{16}{5} = 3.2$} \\
\hline

\end{tabular}
\label{tab1}
\end{table}

\end{landscape}


\begin{thebibliography}{99}


\bibitem{Bagina_2011}
O~.Bagina, Tiling the plane with convex pentagons (in Russian), 
\textit{Bulletin of Kemerovo State Univ}, \textbf{4} no. 48 (2011) 63--73.
Available online: \url{http://bulletin.kemsu.ru/Content/documents/Bulletin_Kemsu_11_4.pdf}
(accessed on 7 November 2015).

\bibitem{Gardner_1975a}
M.~Gardner, On tessellating the plane with convex polygon tiles, 
\textit{Scientific American}, \textbf{233} no. 1 (1975) 112--117.

\bibitem{Gardner_1975b}
\textemdash, A random assortment of puzzles, 
\textit{Scientific American}, \textbf{233} no. 6 (1975) 116--119.

\bibitem{G_and_S_1987}
B.~Gr\"{u}nbaum, G.C.~Shephard, \textit{Tilings and Patterns}. 
W. H. Freeman and Company, New York, 1987. 
pp.15--35 (Chapter 1), pp.113--157 (Chapter 3), pp.471--487, 
pp.492--497, pp.517--518 (Chapter 9), pp.531--549, and pp.580--582 (Chapter 10).

\bibitem{Hallard_1991}
T.C.~Hallard, J.F.~Kennet, K.G.~Richard, \textit{Unsolved Problems in Geometry}. 
Springer-Verlag, New York, 1991. pp.79--80, pp.95--96 (C14), and pp.101--103 (C18).

\bibitem{Hirshh_1985}
M.D.~Hirschhorn, D.C.~Hunt, Equilateral convex pentagons which tile the plane, 
\textit{Journal of Combinatorial Theory}, Series A \textbf{39} (1985) 1--18.

\bibitem{Kershner_1968}
R.B.~Kershner, On paving the plane, 
\textit{American Mathematical Monthly}, \textbf{75} (1968) 839--844.

\bibitem{Klamkin_1980}
M.S.~Klamkin, A.~Liu, Note on a result of Niven on impossible tessellations, 
\textit{American Mathematical Monthly}, \textbf{87} (1980) 651--653.

\bibitem{Mann_2015}
C.~Mann, J.~McLoud-Mann, D.~Von~Derau, Convex pentagons that 
admit i-block transitive tilings (2015), 
\url{http://arxiv.org/abs/1510.01186} (accessed on 23 October 2015).

\bibitem{Reinhardt_1918}
K.~Reinhardt, \textit{\"{U}ber die Zerlegung der Ebene in Polygone}, 
Inaugural-Dissertation, Univ. Frankfurt a.M., R. Noske, Boran and Leipzig, 1918.
Available online: \url{http://gdz.sub.unigoettingen.de/dms/load/img/?PPN=PPN316479497&IDDOC=57097}
(accessed on 7 November 2015).

\bibitem{Schatt_1978}
D.~Schattschneider, Tiling the plane with congruent pentagons, 
\textit{Mathematics Magazine}, \textbf{51} (1978) 29--44.  

\bibitem{Stein_1985}
R.~Stein, A new pentagon tiler, \textit{Mathematics magazine}, \textbf{58} (1985) 308.

\bibitem{Sugimoto_2012a}
T.~Sugimoto, Convex pentagons for edge-to-edge tiling, I, 
\textit{Forma}, \textbf{27} (2012) 93--103.
Available online: \url{http://www.scipress.org/journals/forma/abstract/2701/27010093.html}
(accessed on 8 August 2015).

\bibitem{Sugimoto_2012b}
\textemdash, Convex pentagons for edge-to-edge tiling and convex polygons 
for aperiodic tiling (in Japanese), \textit{Ouyou suugaku goudou kenkyuu syuukai 
houkokusyuu $($Proceeding of Applied Mathematics Joint Workshop$)$} (2012) 
126--131.

\bibitem{Sugimoto_2014}
\textemdash, Consideration of Periodic Tilings by Convex Pentagons (in Japanese), 
\textit{Katachi no kagaku kaishi $($Bulletin of The Society for Science on Form$)$} (2014) 
48--49.

\bibitem{Sugimoto_2015}
\textemdash, Convex pentagons for edge-to-edge tiling, I\hspace{-.1em}I, 
\textit{Graphs and Combinatorics}, \textbf{31} (2015) 281--298. 
doi: 10.1007/s00373-013-1385-x.
Available online: \url{http://dx.doi.org/10.1007/s00373-013-1385-x}
(accessed on 8 August 2015).

\bibitem{Sugimoto_NoteTP}
\textemdash, Tiling Problem: Convex pentagons for edge-to-edge tiling and 
convex polygons for aperiodic tiling (2015).
\url{http://arxiv.org/abs/1508.01864} (accessed on 16 November 2015).

\bibitem{Sugimoto_2016}
\textemdash, Convex pentagons for edge-to-edge tiling, I\hspace{-.1em}I\hspace{-.1em}I, 
\textit{Graphs and Combinatorics}, \textbf{32} (2016) 785--799.  
doi:10.1007/s00373-015-1599-1.
Available online: \url{http://dx.doi.org/10.1007/s00373-015-1599-1}
(accessed on 8 August 2015).

\bibitem{Sugimoto_APTCP}
\textemdash, Convex polygons for aperiodic tiling (2016).
\url{http://arxiv.org/abs/1602.06372} (accessed on 20 February 2016).

\bibitem{Sugi_Ogawa_2006}
T.~Sugimoto, T.~Ogawa, Properties of tilings by convex pentagons, 
\textit{Forma}, \textbf{21} (2006) 113--128.
Available online: \url{http://www.scipress.org/journals/forma/abstract/2102/21020113.html}
(accessed on 8 August 2015).




\end{thebibliography}
\end{document}